\documentclass[final,1p,times]{elsarticle}
\usepackage{graphicx}
\usepackage{epsfig}
\usepackage{psfrag}
\usepackage{amstext,amsthm,amsfonts,amsmath,amsthm,graphicx,amssymb,amscd,amsbsy}
\newtheorem{thm}{Theorem}
\newtheorem{cor}{Corollary}
\newtheorem{lem}{Lemma}
\newtheorem{prop}{Proposition}
\newdefinition{defn}{Definition}
\newdefinition{rmk}{Remark}

\numberwithin{equation}{section}

 \usepackage{lineno}

\begin{document}

\begin{frontmatter}
\title{Asymptotic stability at infinity for bidimensional\\ Hurwitz vector fields \tnoteref{cor1}}
\author{Roland Rabanal\fnref{R. Rabanal}}
\tnotetext[cor1]{This paper was written when the author served as an Associate Fellow at \textsc{ictp-Italy}.}
\fntext[{R. Rabanal}]{The author was partially supported by \textsc{pucp-Peru}~(\textsc{dai}: 2012-0020).}
\ead{rrabanal@pucp.edu.pe}
\begin{abstract}
Let $X:U\to\mathbb{R}^2$  be a differentiable vector field. Set $\mbox{\rm Spc}(X)=\{\mbox{eigenvalues of } DX(z) : z\in U\}$.
This $X$ is called Hurwitz if $\mbox{\rm Spc}(X)\subset\{z\in\mathbb{C}:\Re(z)<0\}$.
Suppose that $X$ is Hurwitz and $U\subset\mathbb{R}^2$ is the complement of a compact set. Then by adding to $X$ a constant $v$ one obtains that the infinity is either an attractor or a repellor for $X+v.$ That means:
\mbox{(i) there} exists an unbounded sequence of closed curves, pairwise bounding an annulus the boundary of which  is transversal to $X+v$, and
\mbox{(ii) there} is a neighborhood of infinity with unbounded trajectories, free of singularities and periodic trajectories of $X+v.$
This result is obtained after to proving  the existence of $\tilde{X}:\mathbb{R}^2\to\mathbb{R}^2,$ a topological embedding such that $\tilde{X}$ equals $X$ in the complement of some compact subset of $U$.
\end{abstract}


\end{frontmatter}

\linenumbers
\section{Introducion}\label{sec:introd}
A basic example of non--discrete dynamics on the Euclidean space is given by a linear vector field.
This linear system is infinitesimally hyperbolic if every eigenvalue has nonzero real part, and it has well known properties \cite{Har02,Chicone06}. For instance, when the real part of its eigenvalues are negative (Hurwitz matrix), the origin is a global attractor rest point.
In the nonlinear case, there has been a great interest in the local study of vector fields around their rest points~\cite{Du,Rou,DSR,Ta74}.
However, in order to describe a global phase-portrait, as in \cite{piresrabanal13,R09center,Fesler,glutsyuk,Gutierrez95} it is absolutely necessary to study its behavior in a neighborhood of infinity \cite{Olech-3}.
\par
The Asymptotic Stability at Infinity has been investigated with a strong influence of \cite{Olech-1}, where Olech showed a connection between stability and injectivity (see also \cite{GJLT07,FGR1,GV,R10,piresrabanal13}).
This research was also studied in \cite{Olech-3,GT95,GS03,GPR06,R06,AGG06}.
In \cite{GT95}, Gutierrez and Teixeira study $C^1-$vector fields $Y:\mathbb{R}^2\to\mathbb{R}^2$, the linearizations of which satisfy
\mbox{(i) $\mbox{\rm det}(DY(z))>0$} and \mbox{(ii) $\mbox{\rm Trace}(DY(z))<0$} in an neighborhood of infinity.
By using \cite{Gutierrez95}, they prove that if $Y$ has a rest point
and the Index $\mathcal{I}(Y)=\int\mbox{\rm Trace}~(DY)<0$ (resp. $\mathcal{I}(Y)\geq 0$), then $Y$ is topologically equivalent to $z\mapsto -z$ that is \lq\lq the infinity is a repellor \rq\rq (resp. to $z\mapsto z$ that is \lq\lq the infinity is an attractor\rq\rq).
This Gutierrez-Teixeira's paper was used to obtain the next theorem, where $\mbox{\rm Spc}(Y)=\{\mbox{eigenvalues of } DY (z):z\in \mathbb{R}^2\setminus\overline{D}_\sigma\}$, and $\Re(z)$ is the real part of $z\in\mathbb{C}$.

\begin{thm}[Gutierrez-Sarmiento]\label{thm:GS03}
Let  $Y:\mathbb{R}^2\setminus\overline{D}_{\sigma}\to\mathbb{R}^2$
be a $C^1-$ map, where $\sigma>0$ and
$\overline{D}_\sigma=\{z\in\mathbb{R}^2:||z||\leq\sigma\}.$
The following is satisfied:
\begin{enumerate}
    \item [(i)]If for some $\varepsilon>0,$
    $\mbox{\rm Spc}(Y)\cap(-\varepsilon,+\infty)=\emptyset.$
    Then there exists
    $s\geq\sigma$ such that the restriction  $Y|:{{\mathbb{R}^2\setminus \overline{D}_s}}\to \mathbb{R}^2$ is
    injective.
    \item [(ii)]If for some $\varepsilon>0,$ the spectrum $\mbox{\rm Spc}(Y)$ is disjoint of
    the union
    $(-\varepsilon,0]\cup\{z\in\mathbb{C}:\Re(z)\geq0\}$.
    Then there exist $p_0\in\mathbb{R}^2$ such that the point
    $\infty$ of the Riemann Sphere $\mathbb{R}^2\cup\{\infty\}$ is
    either an attractor or a repellor of $z^\prime=Y(z)+p_0.$
\end{enumerate}
\end{thm}

\mbox{Theorem \ref{thm:GS03}} is given in \cite{GS03}, and it has been extended to  differentiable maps $X:\mathbb{R}^2\setminus\overline{D}_\sigma\to\mathbb{R}^2$ in \cite{GR06,GPR06}.
In both papers  the eigenvalues also avoid a real open neighborhood of zero.
In  \cite{R06} the author examine the intrinsic relation between the asymptotic behavior of $\mbox{\rm Spc}(X)$ and the global injectivity of the local diffeomorphism given by $X$.
He uses $Y_{\theta}=R_{\theta}\circ{Y}\circ R_{-\theta}$, where $R_{\theta}$ is the linear rotation of angle $\theta\in\mathbb{R}$, and (motivated by \cite{GV}) introduces the so--called $B-$condition \cite{R06e,R10}, which claims:
\begin{description}
\item[] \emph{for each $\theta\in\mathbb{R}$, there does not exist a sequence $(x_k,y_k)\in\mathbb{R}^2$ with $x_k\to+\infty$ such that ${Y_\theta}((x_k,y_k))\rightarrow p\in\mathbb{R}^2$ and $DY_{\theta}(x_k,y_k)$ has a real eigenvalue $\lambda_k$ satisfying $x_k\lambda_k\rightarrow 0$.}
\end{description}
By using this, \cite{R06} improves the differentiable version of Theorem~\ref{thm:GS03}.
\par
In the present paper we prove that the condition
\[\mbox{Spc}(X)\subset\{z\in\mathbb{C}:\Re(z)<0\}\]
is enough in order to obtain Theorem~\ref{thm:GS03} for differentiable vector fields $X:\mathbb{R}^2\setminus\overline{D}_\sigma\to\mathbb{R}^2$.
\par
Throughout this paper, $\mathbb{R}^2$ is embedded in the Riemann sphere $\mathbb{R}^2\cup\{\infty\}$.
Thus $(\mathbb{R}^2\setminus\overline{D}_\sigma)\cup\{\infty\}$ is the subspace of  $\mathbb{R}^2\cup\{\infty\}$ with the induced topology, and \lq infinity \rq  refers to the point $\infty$ of $\mathbb{R}^2\cup\{\infty\}$.
Moreover given $C\subset\mathbb{R}^2$, a closed (compact, no boundary) curve ($1-$manifold), $\overline{D}(C)$ (respectively $D(C)$) is the compact disc (respectively open disc) bounded by $C$.
Thus, the boundaries $\partial \overline{D}(C)$ and $\partial D(C)$ are equal to $C$, homeomorphic to $\partial D_1=\{z\in\mathbb{R}^2:||z||=1\}$.

\section{Statements of the results}\label{sec:2}
For every $\sigma>0$ let $\overline{D}_\sigma=\{z\in\mathbb{R}^2:||z||\leq\sigma\}$.
Outside this compact disk we consider a differentiable vector field
$X:\mathbb{R}^2\setminus\overline{D}_{\sigma}\to\mathbb{R}^2$.
As usual, a \textbf{trajectory of $X$ starting at $q\in\mathbb{R}^2\setminus\overline{D}_{\sigma}$} is defined as the integral curve determined by a maximal solution of the initial value problem $\dot{z}=X(z),~ z(0)=q$.
This is a curve $I_q\ni t\mapsto\gamma_q(t)=(x(t),y(t))$, satisfying:
\begin{itemize}
  \item $t$ varies on some open real interval containing the zero, the image of which $\gamma_q(0)=q$;
  \item $\gamma_q(t)\in\mathbb{R}^2\setminus\overline{D}_{\sigma}$ and there exist the real derivatives $\frac{dx}{dt}(t)$, $\frac{dy}{dt}(t)$;
  \item $\dot{\gamma_q}(t)=\Big(\frac{dx}{dt}(t),\frac{dy}{dt}(t)\Big)$ the velocity vector field of $\gamma_q$ at  $\gamma_q(t)$ equals $X(\gamma_q(t))$ and
  \item  $I_q\subset\mathbb{R}$ is the maximal interval of definition.
\end{itemize}
We identify the trajectory $\gamma_q$ with its image $\gamma_q(I_q)$, and we denote by ${\gamma^+_q}$ (resp. ${\gamma^-_q}$) the positive (resp. negative) semi-trajectory of $X$, contained in $\gamma_q$ and  starting at $q$. In this way $\gamma_q={\gamma_q^-}\cup{\gamma_q^+}$.
Thus  each trajectory has its two limit sets, $\alpha(\gamma_q^-)$ and $\omega(\gamma_q^+)$ respectively.
These limit sets are well defined in the sense that they only depend on the respective solution.
\par
A $C^0-$vector field $X:\mathbb{R}^2\setminus\overline{D}_\sigma\to\mathbb{R}^2\setminus\{0\}$ (without rest points) can be extended to a map  \[\hat{X}:(({\mathbb{R}^2\setminus\overline{D}_\sigma}\cup\infty),\infty)\longrightarrow(\mathbb{R}^2,0)\] (which takes $\infty$ to $0$) \cite{AGG06}.
In this manner, all questions concerning the local theory of isolated rest points of $X$ can be formulated and examined in the case of the vector field $\hat{X}$.
For instance, if $\gamma_p^+$ (resp. $\gamma_p^-$) is an unbounded semi-trajectory of $X:\mathbb{R}^2\setminus\overline{D}_\sigma\to\mathbb{R}^2$ starting at $p\in\mathbb{R}^2\setminus\overline{D}_\sigma$ with empty $\omega-$limite (resp. $\alpha-$limit) set, we say ${\gamma_p^+}$ \textbf{goes to infinity} (resp. $\gamma_p^-$ \textbf{comes from infinity}), and it is denoted by $\omega({\gamma_p^+})=\infty$ (resp. $\alpha({\gamma_p^-})=\infty$).
Therefore, we may also talk about the phase portrait of $X$ in a neighborhood of $\infty$.
\par
As in our paper \cite{GPR06}, we say  that the point  at infinity $\infty$ of the Riemann Sphere $\mathbb{R}^2\cup\{\infty\}$ is an \textbf{attractor} (resp. a \textbf{repellor}) for the continuos vector field $X:\mathbb{R}^2\setminus\overline{D}_\sigma\to\mathbb{R}^2$ if:
\begin{itemize}
    \item There is a sequence of closed curves, transversal  to $X$ and tending to infinity.
        That is for every $r\geq\sigma$ there exist a closed curve $C_r$ such that $D(C_r)$ contains $D_r$ and $C_r$ has transversal contact to each small local integral curve of $X$ at any $p\in C_r$.
    \item For some $C_s$ with $s\geq\sigma,$ all
    trajectories $\gamma_p$ starting at a point $p\in
    \mathbb{R}^2\setminus\overline{D}(C_s)$ satisfy $\omega({\gamma^+_p})=\infty$ that is $\gamma_p^+$ go to infinity (resp. $\alpha({\gamma^-_p})=\infty$ that is  $\gamma_p^-$ come from infinity).
\end{itemize}
We also recall that $\mathcal{I}(X)$, the \textbf{index of $X$ at infinity} is the number of the extended line $[-\infty,+\infty]$ given by
\[
\mathcal{I}(X)={\int_{\mathbb{R}^2}}\mbox{\rm
Trace}(D\widehat{X})dx\wedge dy,
\]
where $\widehat{X}:\mathbb{R}^2\to\mathbb{R}^2$ is a global differentiable vector field such that:
\begin{itemize}
    \item In the complement of some disk $\overline{D}_s$ with $s\geq\sigma$ both $X$ and $\widehat{X}$ coincide.
    \item The map $z\mapsto\mbox{\rm Trace}(D\widehat{X}_z)$ is Lebesgue almost--integrable in whole $\mathbb{R}^2,$ in the sense of \cite{GPR06}.
\end{itemize}
This index is a well-defined number in $[-\infty,+\infty]$, and it does not depend on the pair  $(\widehat{X},s)$ as shown \cite[\mbox{Lemma 12}]{GPR06}.
\begin{defn}
The differentiable vector field (or map) $X:\mathbb{R}^2\setminus{\overline{D}_\sigma}\to\mathbb{R}^2$ is called \textbf{Hurwitz} if every eigenvalue of the Jacobian matrices has negative real part.
This means that its spectrum  $\mbox{\rm Spc}(X)=\{\mbox{eigenvalues of } DX (z):z\in \mathbb{R}^2\setminus\overline{D}_\sigma\}$ satisfies $\mbox{\rm Spc}(X)\subset\{z\in\mathbb{C}:\Re(z)<0\}$.
\end{defn}

We are now ready to state our result.

\begin{thm}\label{thm:main-1}
Let  $X:\mathbb{R}^2\setminus{\overline{D}_\sigma}\to\mathbb{R}^2$ be  a differentiable  vector field (or map),  where $\sigma>0$ and ${\overline{D}_\sigma}=\{z\in\mathbb{R}^2:||{z}||\leq\sigma\}$.
Suppose that $X$ is Hurwitz: $\mbox{\rm  Spc}(X)\subset\{z\in\mathbb{C}:\Re(z)<0\}$. Then
\begin{enumerate}
\item[(i)] There are $s\geq\sigma$ and a globally injective local homeomorphism $\widetilde{X}:\mathbb{R}^2\to\mathbb{R}^2$ such that $\widetilde{X}$ and $X$ coincide on $\mathbb{R}^2\setminus{\overline{D}_s}$.
    Moreover, the restriction  $X|:{{\mathbb{R}^2\setminus \overline{D}_s}}\to \mathbb{R}^2$ is injective, and it admits a global differentiable extension $\widehat{X}$ such that the pair $(\widehat{X},s)$ satisfies the definition of the index of $X$ at infinity, and this index $\mathcal{I}(X)\in[-\infty,+\infty)$.
\item[(ii)] For all $p\in\mathbb{R}^2\setminus{\overline{D}_\sigma}$, there is a unique positive semi-trajectory of $X$ starting  at $p$.
    Moreover, for some $v\in\mathbb{R}^2$, the point at infinity $\infty$ of the Riemann Sphere $\mathbb{R}^2\cup\{\infty\}$ is an attractor (respectively a repellor) for the vector field $X+v:\mathbb{R}^2\setminus\overline{D}_s\to\mathbb{R}^2$ as long as the well-defined index $\mathcal{I}(X)\geq0$ (respectively $\mathcal{I}(X)<0$).
\end{enumerate}
\end{thm}

The map $\widetilde X$ of \mbox{Theorem \ref{thm:main-1}} is not necessarily a homeomorphism.
This $\widetilde X$ is a topological embedding, the image of which may be properly contained in $\mathbb{R}^2$.
Furthermore, if $A:\mathbb{R}^2\to\mathbb{R}^2$ is an arbitrary invertible linear map, then \mbox{Theorem \ref{thm:main-1}} applies to the map $A\circ{X}\circ{A^{-1}}$.
\par
\mbox{Theorem \ref{thm:main-1}} improves the main results of \cite{GR06,GPR06}.
\mbox{Item (i)} complements the injectivity work of \cite{GR06} (see also \cite{R06,GS03}), where the authors consider the assumption $\mbox{\rm Spc}(X)\cap(-\varepsilon,+\infty)=\emptyset$ (as in \mbox{Theorem \ref{thm:GS03}}).
\mbox{Item (ii)} generalizes \cite{GPR06}, where the authors utilize the second condition of \mbox{Theorem \ref{thm:GS03}}.
In our new assumptions, the negative eigenvalues can tend to zero.

\subsection{Description of the proof of \mbox{Theorem \ref{thm:main-1}}}
Since the Local Inverse Function Theorem is true, a map $X=(f,g):\mathbb{R}^2\setminus\overline{D}_{\sigma}\to\mathbb{R}^2$ as in \mbox{Theorem \ref{thm:main-1}} is a local diffeomorphism.
Thus the level curves $\{f={\rm constant}\}$ make up a  $C^0-$foliation $\mathcal{F}(f)$ the leaves of which are differentiable curves, and the restriction of the other submersion $g$ to each of these leaves is strictly monotone.
In particular, $\mathcal{F}(f)$ and $\mathcal{F}(g)$ are (topologically) transversal to each other.
We orient $\mathcal{F}(f)$  in agreement that if $L_p(f)$ is an oriented leaf of $\mathcal{F}(f)$ thought the point $p$, then the restriction $g\vert:{L_p(f)}\to \mathbb{R}$ is an increasing function in conformity with the orientation of $L_p(f)$.
We denote by $L_p^+(f)=\{z\in L_p(f):g(z)\geq g(p)\}$ (resp. $L_p^+(g)=\{z\in L_p(g):f(z)\geq f(p)\}$) and  $L_p^-(f)=\{z\in L_p(f):g(z)\leq g(p)\}$ (resp. $L_p^-(g)=\{z\in L_p(g):f(z)\leq f(p)\}$) the respective positive and negative half-leaf of $\mathcal{F}(f)$ (resp. $\mathcal{F}(g)$). Thus ${L_q}(f)={L_q^-}(f)\cup{L_q^+}(f)$ and ${L_q^-}(f)\cap{L_q^+}(f)=\{q\}$.
In this context,  the nonsingular vector fields
\begin{equation}\label{eq:hamiltonian}
X_f=(-f_y,f_x)\quad\mbox{\rm and}\quad\widetilde{X}_g=(g_y,-g_x),
\end{equation}
given by the partial derivatives are tangent to $L_p(f)$ and $L_p(g)$, respectively.
This construction has previously been used in \cite{FGR1}.
\par
We need the following defini\-tion \cite{GS03,FGR1}. Let $h_0(x,y)=xy$ and consider the set
\[B=\Big\{(x,y)\in [0,2]\times[0,2]:0<x+y\leq 2 \Big\}.\]
\begin{defn}\label{def:HRC}
Let $X=(f,g)$ be a differentiable local homeomorphism.
Given $h\in\{f,g\}$, we say that $\mathcal{A}$ (in the domain of $X$) is a \textbf{half-Reeb component} for $\mathcal{F}(h)$ if there is a homeomorphism $H:B \to \mathcal{A}$ which is a topological equivalence between $\mathcal{F}(h) \vert_{\mathcal{A}}$ and ${\mathcal{F}}(h_0) \vert_B$ such that: (\mbox{Fig. \ref{fig:HRChurwitz}})
\begin{itemize}
    \item The segment $\{(x,y)\in B : x+y=2\}$ is sent by $H$ onto a transversal section for the foliation $\mathcal{F}(h)$ in the complement of the point $H(1,1)$; this section is called the {\rm compact edge} of $\mathcal{A}$.
    \item Both segments $\{(x,y)\in B : x=0 \}$ and $\{(x,y)\in B : y=0\}$ are sent by $H$ onto full half-leaves of $\mathcal{F}(h)$.
        These half-leaves of $\mathcal{F}(h)$ are called the {\rm non--compact edges} of $\mathcal{A}$.
\end{itemize}
Observe that $\mathcal{A}$ may not be a closed subset of $\mathbb{R}^2$, and $H$ does not need to be extended to infinity.
\end{defn}
\par
\mbox{Section \ref{sec:3}} gives new results on the foliations induced by a local diffeomorphism $X=(f,g):\mathbb{R}^2\setminus\overline{D}_{\sigma}\to\mathbb{R}^2$ (see also Proposition~\ref{prop:10}). \mbox{Theorem \ref{teo:1}} implies that the conditions
\begin{equation}\label{eq:cond}
\mbox{\rm Spc}(X)\cap[0,+\infty)=\emptyset \mbox{ and \lq each half--Reeb component of } \mathcal{F}(f) \mbox{ is bounded \rq}
\end{equation}
give the existence of $s\geq\sigma$ such that $X|_{\mathbb{R}^2\setminus{\overline{D}_s}}$ can be extended to an injective map $\widetilde{X}:\mathbb{R}^2\to\mathbb{R}^2$.
\mbox{Section \ref{sec:4}} presents some preliminary results on maps such that $\mbox{\rm Spc}(X)\cap[0,+\infty)=\emptyset$.
Section~\ref{sec:5} concludes with the proof of \mbox{Theorem \ref{thm:main-1}}. The main step is given in \mbox{Proposition \ref{prop:main}}, which implies that Hurwitz maps satisfy \eqref{eq:cond}.
Therefore, \mbox{Theorem \ref{thm:main-1}} is obtained by using this \mbox{Proposition \ref{prop:main}} and some previous work \cite{GR06,GPR06}.

\section{Local diffeomorphisms that are injective on unbounded open sets}\label{sec:3}
Let $X=(f,g):\mathbb{R}^2\setminus\overline{D}_{\sigma}\to\mathbb{R}^2$ be an orientation preserving local diffeomorphism, that is $\mbox{\rm det}(DX)>0$.
Next subsection gives preparatory results about $\mathcal{F}(f)$ in order to obtain that $X$ will be injective on topological half planes (see Proposition~\ref{prop:10}).
Subsection~\ref{subsec:2.3} presents a condition under which $X$  is injective at infinity, that is outside some compact set.

\subsection{Avoiding tangent points}
Let $C\subset\mathbb{R}^2\setminus{\overline{D}_\sigma}$ be a closed curve surrounding the origin.
We say that the vector field $X:\mathbb{R}^2\setminus{\overline{D}_\sigma}\to\mathbb{R}^2$ has \textit{contact} (resp. tangency with; resp. transversal to; etc) with $C$ at $p\in C$ if for each small local integral curve of $X$ at $p$ has such property.

\begin{defn}\label{def:general-position}
A closed curve $C\subset\mathbb{R}^2\setminus{\overline{D}_\sigma}$ is in \textbf{general--position} with $\mathcal{F}(f)$ if there exists a set $T\subset C$, at most finite such that:
\begin{itemize}
  \item $\mathcal{F}(f)$ is transversal to $C\setminus T$.
  \item $\mathcal{F}(f)$ has a tangency with $C$ at every point of $T$.
  \item  A  leaf of $\mathcal{F}(f)$ can meet tangentially $C$ at most at one point.
\end{itemize}
\end{defn}

Denote by $\mathcal{GP}(f,s)$ the set of all closed curves $C\subset\mathbb{R}^2\setminus{\overline{D}_\sigma}$ in general--position with $\mathcal{F}(f)$  such that $D_s\subset D(C)$.
If $C\subset\mathbb{R}^2\setminus{\overline{D}_\sigma}$ is in general--position, we denote by  $n^e(C,f)$ (resp. $n^i(C,f)$) the number of tangent points of $\mathcal{F}(f)$ with $C$, which are external (resp. internal).
Here, external (resp. internal) means the existence of a small open interval $(\tilde{p},\tilde{q})_f\subset L_p(f)$ such that the intersection set $(\tilde{p},\tilde{q})_f\cap C=\{\mbox{tangent point}\}$ and $(\tilde{p},\tilde{q})_f\subset\mathbb{R}^2\setminus D(C)$ (resp. $(\tilde{p},\tilde{q})_f\subset\overline{D}(C)$).

\begin{rmk}\label{rem:a,b}
If $C\subset\mathbb{R}^2\setminus{\overline{D}_\sigma}$ is in general--position with $\mathcal{F}(f),$ there exist $a,b\in C$ two different point such that $f(C)=[f(a),f(b)].$
Moreover, $a$ and $b$ are external tangent points because the map $(f,g)$ is an orientation preserving local diffeomorphism.
Since $C$ and $f(C)$ are connected and $C$ is not contained in any leaf  of $\mathcal{F}(f)$, we conclude that both external tangencies are different.
\end{rmk}

Corollary~\ref{cor:ext}  gives important properties of the leafs passing trough a point in Remark~\ref{rem:a,b}, if we select $C\in\mathcal{GP}(f,\sigma)$ with the minimal number of internal tangencies.  To this end, the next lemma will be needed.

\begin{lem}\label{lem:external-tangency}
Let $C\in\mathcal{GP}(f,\sigma)$.
Suppose that a leaf  $L_q(f)$ of $\mathcal{F}(f)$ meets $C$ transversally somewhere and with an external tangency at a point ${p\in C}$.
Then $L_q(f)$ contains a closed subinterval $[p,r]_f$ which meets $C$ exactly at $\{p,r\}$ (doing it transversally at $r$) and the following is satisfied:
\begin{itemize}
\item[(a)] If $[p,r]$ is the closed subinterval of $C$ such that $\Gamma=[p,r]\cup[p,r]_f$ bounds a compact disc $\overline{D}(\Gamma)$ contained in $\mathbb{R}^2\setminus D(C),$ then points of $L_q(f)\setminus[p,r]_f$ nearby $p$ do not belong to $\overline{D}(\Gamma)$.
\item[(b)] Let $(\tilde{p},\tilde{r})$ and $[\tilde{p},\tilde{r}]$ be subintervals of  $C$ satisfying $[p,r]\subset(\tilde{p},\tilde{r}) \subset [\tilde{p},\tilde{r}]$.
        If $\tilde{p}$ and $\tilde{r}$ are close enough to $p$ and $r,$ respectively; then we may deform $C$ into $C_1\in\mathcal{GP}(f,\sigma)$ in such a way that the deformation fixes $C\setminus(\tilde{p},\tilde{r})$ and takes $[\tilde{p},\tilde{r}]\subset C$ to a closed subinterval $[\tilde{p},\tilde{r}]_1 \subset C_1$ which is close to $[p,r]_f$.
    Furthermore, the number of generic tangencies  of $\mathcal{F}(f)$ with $C_1$ is smaller than that of $\mathcal{F}(f)$ with $C$.
\end{itemize}
\end{lem}

\begin{proof}
 We refer the reader to  \cite[Lemma~2]{GS03}.
\end{proof}

\begin{cor}\label{cor:ext}
Let $C\in\mathcal{GP}(f,s).$ Suppose that $C$ minimizes ${n^i}(C,f)$ and $f(C)=[f(a),f(b)],$ then $C\cap L_a(f)=\{a\}$ and $C\cap L_b(f)=\{b\}.$
\end{cor}

\begin{proof}
We only consider the case of the point $a$.
Assume by contradiction that the number of elements in $C\cap L_a(f)$ is greater than one i.e $\sharp(C\cap L_a(f))\geq 2$.
The last condition  of \mbox{Definition \ref{def:general-position}} implies that the intersection of $L_a(f)$ with the other point is transversal to $C$.
Then there is a disk as in statement (a) of \mbox{Lemma \ref{lem:external-tangency}}.
By using the second part of  \mbox{Lemma \ref{lem:external-tangency}} we can avoid the external tangency $a$ and some internal tangency.
This is a contradiction because ${n^i}(C,f)$ is minimal.
\end{proof}

\begin{rmk}
Corollary~\ref{cor:ext} remains true if we take any external tangency (not necessarily $a$ and $b$) in a closed curve $C_s\in\mathcal{GP}(f,s)$ with minimal $n^i(C_s,f)$.
\end{rmk}

\subsection{Minimal number of internal tangent points}
A oriented leaf of $\mathcal{F}(f)$ whose distance to $\overline{D}_\sigma$ is different from zero has unbounded half--leaves.
Given any $L_p(f)$ with unbounded half--leaves, we denote by $H^+(L_p)$  and $H^-(L_p)$ the two components of $\mathbb{R}^2\setminus L_p(f)$ in order that $L_p^+(g)\setminus\{p\}$ be contained in $H^+(L_p).$
Therefore, the image $X(H^+(L_p))$ is an open connected subset of the semi--plane $\{x>f(p)\}:=\{(x,y)\in\mathbb{R}^2:x>f(p)\}.$

\begin{rmk}\label{rem:1}
If the image $X(H^+(L_p))$ is a vertical convex set,  all the level curves $\{f=c\}\subset H^+(L_p)$ are connected.
Thus the restriction $X|:H^+(L_p)\to X(H^+(L_p))$ is an homeomorphism, and it sends every leaf of $\mathcal{F}(f){|_{H^+(L_p)}}$ over vertical lines. Therefore, it is a topological equivalence between two foliations.
\end{rmk}

\begin{lem}\label{lem:non-convex.hrc}
If $H^+(L_p)$ is disjoint from $\overline{D}_{\sigma}$ and the image $X(H^+(L_p))$ is not a vertical convex set, then $H^+(L_p)$ contains a  half-Reeb component of $\mathcal{F}(f).$
\end{lem}

\begin{proof}
By remark~\ref{rem:1}, some level set $\{f=\tilde{c}\}\subset H^+(L_p)$  is disconnected. Therefore, the result is obtained directly from \cite[Proposition~1]{GV}.
\end{proof}
Lemma~\ref{lem:non-convex.hrc} and Remark~\ref{rem:1} hold when we
consider $H^-(L_p).$
\begin{lem}\label{lem:def-eta}
Recall that $\mathcal{GP}(f,s)$ is the set of all closed curves $C\subset\mathbb{R}^2\setminus{\overline{D}_\sigma}$ in general--position with $\mathcal{F}(f)$  such that $\overline{D}_s\subset D(C).$
Let ${\eta^i}:[\sigma,\infty)\to\mathbb{N}\cup\{0\}$ be the function given by ${\eta^i}(s)={n^i}(C_s,f)$ where $C_s\in \mathcal{GP}(f,s)$ minimizes the number of internal tangent points with $\mathcal{F}(f).$
The following statements hold:
\begin{enumerate}
\item[(a)] The function $\eta^i$ is nondecreasing.
\item[(b)] If $\eta^i$ is bounded then, there exist $s_0\in[\sigma,\infty)$ such that ${\eta^i}(s)\leq {\eta^i}(s_0)$ for all
    $s\in[\sigma,\infty).$
\item[(c)] Set $f(C_{s_0})=[f(a),f(b)]$. Suppose that
    $\mathcal{F}(f)$ has a half-Reeb component $\mathcal{A}$ whose
    image $f(\mathcal{A})$ is disjoint from $(f(a),f(b)),$ then such
    $s_0\in[\sigma,\infty)$ is not a maximum value of the function
    $\eta^i.$
\end{enumerate}
\end{lem}

\begin{proof}
As $C_{s+1}$ also belongs to $\mathcal{GP}(f,s)$ we have that $\eta^i(s)\leq n^i(C_{s+1},f)$.
Therefore (a) is true.
\par
To prove the second part, we introduce the set $A_s=\{n\in\mathbb{N}\cup\{0\}:n\geq{n^i}(C_s,f)=\eta^i(s)\}$, for every $s\geq\sigma$.
From this definition it is not difficult to check that: $\eta^i$ is bounded if and only if $\displaystyle{\cap_{s\geq\sigma}}A_s\neq\emptyset$.
Therefore, the first element of $\displaystyle{\cap_{s\geq\sigma}}A_s\neq\emptyset$ is the bound $\eta^i(s_0)$ of statement (b).
This proves the second statement.
\par
We shall have established (c) if we prove that there is some  $C_{{s_0}+\varepsilon}$ with $\varepsilon>0$ such that $\eta^i(s_0)+1\leq {n^i}(C_{{s_0}+\varepsilon},f)$.
To this end, we select $s>s_0$ large enough for which $C_s$ is enclosing $\overline{D}(C_{s_0})\cup\Gamma$ where $\Gamma$ is the compact edge of $\mathcal{A}$.
Since, this $C_s$ intersects both leaves $L_p(f)$ and $L_q(f)$ where $p$ and $q$ are the endpoints of $\Gamma$, we obtain that $n^i(C_s,f)$ is greater than $n^i(C_{s_0},f)=\eta^i(s_0)$.
Therefore, for some  $\varepsilon>0$ there is $C_{{s_0}+\varepsilon}$ such that $\eta^i(s_0)+1\leq {n^i}(C_{{s_0}+\varepsilon},f)$. This proves (c).
\end{proof}

\begin{prop}\label{prop:10}
Let
$X=(f,g):\mathbb{R}^2\setminus\overline{D}_{\sigma}\to\mathbb{R}^2$
be a map with  $\mbox{\rm det}(DX)>0.$ Consider  $C_s$ and
$\eta^i$ as in Lemma~\ref{lem:def-eta}. If $C_{s_0}$ satisfies
that ${n^i}(C_s,f)\leq{n^i}({C_{s_0}},f)$ for all
$s\in[\sigma,\infty)$ and $f(C_{s_0})=[f(a),f(b)].$ Then, for each
$p\in\{a,b\}$ at least one of the restrictions of $X$ to
${H^+(L_p)}$ or ${H^-(L_p)}$ is a globally injective map, in
agrement that the domain of this restriction is in the complement
of $\overline{D}(C_{s_0})$.
\end{prop}

\begin{proof}We only consider the case $p=a.$ Suppose that
$H^+(L_a)$ is contained in the complement of
$\overline{D}(C_{s_0}).$ From Remark~\ref{rem:1} it is sufficient
to prove that $X(H^+(L_a))$ is vertical convex. Suppose by
contradiction that it is false. Then
Lemma~\ref{lem:non-convex.hrc} implies that there is a half-Reeb
component $\mathcal{A}\subset H^+(L_a).$  By using statement (c)
of Lemma~\ref{lem:def-eta} this  $s_0$ is not a maximum value of
the function $\eta^i$. This contradiction with our selection of
the circle $C_{s_0}$ conclude the proof.
\end{proof}
\subsection{Extending maps to topological embeddings}\label{subsec:2.3}

The next theorem implies the injectivity  at infinity of a map, and it is obtained by using the methods, ideas and arguments of \cite{GR06}.
We only give the proof, in the case of continuously differentiable maps.

\begin{thm}\label{teo:1}
Let $X=(f,g):\mathbb{R}^2\setminus{\overline{D}_\sigma}\to\mathbb{R}^2$ be an differentiable local homeomorphism with $\mbox{det}(DX)>0$. Suppose that $\mbox{\rm Spc}(X)\cap[0,+\infty)=\emptyset$, and each half--Reeb component of either $\mathcal{F}(f)$ or $\mathcal{F}(g)$ is bounded.
Then there exist $s\geq\sigma$ such that the restriction $X|:\mathbb{R}^2\setminus{\overline{D}_s}\to\mathbb{R}^2$ can be extended to a globally injective local homeomorphism $\widetilde{X}=(\tilde{f},\tilde{g}):\mathbb{R}^2\to\mathbb{R}^2$.
\end{thm}

\begin{proof}
We can apply the results of \cite[pp. 166-174]{Har02} to the continuous vector field $X_f$ and obtain that for each closed curve $C\in\mathcal{GP}(f,\sigma)$ the \emph{Index of $X_f$ along $C$,} denoted by $Ind(X_f;C),$ satisfies
\[
Ind(X_f;C)=\frac{2-{{n^e}(f,C)}+{{n^i}(f,C)}}{2}.
\]
If $X_f$ is discontinuous, we proceed as in  \cite{GR06} by using the index of the foliation $\mathcal{F}(f)$ which also satisfies this formulae.
\begin{enumerate}
    \item [(a.1)] We claim that $n^e(f,C)={n^i}(f,C)+2,$ for all $C\in \mathcal{GP}(f,\sigma).$
\end{enumerate}

Suppose that (a.1) is false, so there is $\tilde{C}^1\in \mathcal{GP}(f,\sigma)$  whose
$Ind(X_f;\tilde{C}^1)\neq 0.$ Thus, for some point in $\tilde{C}^1$ the Hamiltonian $X_f(p)=(-f_y(p),f_x(p))$ is vertical.
More precisely we can obtain $p\in \tilde{C}^1$ such that $f_y(p)=0$ and $f_x(p)>0.$  This is a contradiction with the eigenvalue assumptions because $f_x(p)\in \mbox{\rm Spc}(X)\cap (0,+\infty)$.
(If $X_f$ is discontinuous, we refer the reader to \cite[\mbox{Proposition 3.1}]{GR06} where proves that the index of the foliation $\mathcal{F}(f)$ is zero). Therefore, (a.1) holds.

\begin{enumerate}
\item [(a.2)] We claim that if $C_{\sigma}\subset\mathbb{R}^2\setminus\overline{D}_s$
    \emph{ minimizes} ${n^i}(f,C_{\sigma}),$ then every internal tangency in $C_{\sigma}$ gives a half--Reeb component.
\end{enumerate}

For every internal tangency $q\in C_{\sigma}$ we consider the forward Poincar\'{e} map $T:[p,q)_{\sigma}\subset C_{\sigma}\to C_{\sigma}$ induced by the oriended $\mathcal{F}(f)$ (if $T:(q,r]_{\sigma}\subset C_{\sigma}\to C_{\sigma}$ the proof is similar) where $[p,q)_{\sigma}\subset C_{\sigma}$ is the maximal
connected domain of definition of $T$ on which this first return map is continuous.
    If the open arc $L^+_p(f)\setminus\{p\}$ intersects $C_{\sigma}$ we apply Lemma~\ref{lem:external-tangency},
    so we can deform $C_{\sigma}$ in a new circle $\tilde{C}^1\subset\mathbb{R}^2\setminus\overline{D}_s$
    such that the number of internal tangencies of $\tilde{C}^1$ with $\mathcal{F}(f)$ is (strictly) smaller
    than that of $C_{\sigma}.$ This is a contradiction.
    Therefore $L^+_p(f)\setminus\{p\}$ is disjoint from $C_{\sigma}.$
    By using this and our selection of $[p,q)_{\sigma}\subset C_{\sigma}$ is not difficult to check that
    there is a half--Reeb component of $\mathcal{F}(f)$ whose compact
    edge is contained in $C_s,$ Thus, we obtain (a.2).

Notice that, for every circle as in (a.2) any internal tangency of this $C_\sigma$
gives an unbounded half--Reeb component,  thus by our
assumptions and (a.1) we have that  ${n^i}(f,C_\sigma)=0.$ Therefore,
\begin{enumerate}
    \item [(a.3)]
    if $C_{\sigma}\in \mathcal{GP}(f,\sigma)$ is as in (a.2) then ${n^i}(f,C_{\sigma})=0$ and $n^e(f,C_{\sigma})=2.$
    Moreover, $\mathcal{F}(f),$ restricted to $\mathbb{R}^2\setminus D(C_{\sigma}),$ is topologically equivalent
    to the foliation made up by all the vertical straight lines, on $\mathbb{R}^2\setminus D_1$.
\end{enumerate}

Since $X$ has no unbounded half--Reeb component, we can use the last section of \cite{GR06}
(see Proposition 5.1) and obtain that the closed curve $C_{\sigma}$ of (a.3) can be deformed so that,
the resulting new circle $C$ has an exterior collar neighborhood
$U\subset\mathbb{R}^2\setminus D(C)$ such that:
\begin{enumerate}
    \item [(b)] $X(C)$ is a non-trivial closed curve, $X(U)$ is an exterior collar
    neighborhood of $X(C)$ and the restriction $X{|}:U\to X(U)$ is a
    homeomorphism.
\end{enumerate}

By Schoenflies Theorem \cite{B} the map $X{|}:C\to X(C)$ can be extended to a homeomorphism ${X_1}:{\overline{D}(C)}\to{\overline{D}(X(C))}$.
We extend $X:\mathbb{R}^2\setminus D(C)\to \mathbb{R}^2$ to $\widetilde{X}=(\tilde{f},\tilde{g}):\mathbb{R}^2 \to \mathbb{R}^2$ by defining $\widetilde{X}{|_{\overline{D}(C)}}=X_1.$
Thus $\widetilde{X}{|}:U\to X(U)$ is a homeomorphism and $U$ (resp. $X(U)$) is a exterior collar neighborhood of $C$ (resp. $X(C)$).
Consequently, $\widetilde{X}$ is a local homeomorphism and $\mathcal{F}(\tilde{f})$ is topologically equivalent to the foliation made up by all the vertical straight lines.
The injectivity of $\tilde{X}$ follows from the fact that $\mathcal{F}(\tilde{f})$ in trivial \cite[Proposition 1.4]{FGR1}. This concludes the proof.
\end{proof}

\begin{cor}\label{cor:ext with no HRC}
Suppose that $X$ satisfies \mbox{Theorem \ref{teo:1}}.
Then the respective extension  $\tilde{X} = (\tilde f, \tilde g):\mathbb{R}^2\to \mathbb{R}^2$ of $X$ is a globally injective local homeomorphism the foliations of which, $\mathcal{F}(\tilde f)$ and $\mathcal{F}(\tilde g)$ have no half-Reeb components.
\end{cor}
\begin{proof}
We reefer the reader to affirmation (a.3) in the proof of \mbox{Theorem \ref{teo:1}}.
\end{proof}

\begin{cor}
Suppose that $X=(f,g):\mathbb{R}^2\setminus\overline{D}_{\sigma}\to\mathbb{R}^2$ is an orientation preserving local diffeomorphism.
Then the foliation $\mathcal{F}(f)$ (resp. $\mathcal{F}(g)$) has at most countably many half--Reeb components.
\end{cor}

\begin{proof} A half--Reeb component has a tangency with some $C_n$ of Lemma~\ref{lem:def-eta}  with
$n\in\mathbb{N}$ large enough. We conclude, since the closed curves has at most a finite number of tangent points.
\end{proof}

\begin{rmk}
By using a smooth embedding $(f,g):\mathbb{R}^{2}\to\mathbb{R}^{2},$ the authors of \cite[Proposition 1]{GJLT07} prove the existence of foliations  $\mathcal{F}(f)$ which have infinitely many half-Reeb components.
\end{rmk}

\section{Maps free of positive eigenvalues}\label{sec:4}
In this section we present some properties of a map the spectrum of which is disjoint of $[0,+\infty)$.
These results will be used in Section~\ref{sec:5} to proving the first part of Theorem~\ref{thm:main-1}.
In this context, we consider $\mathcal{F}(f)$ and their trajectories $L_q=L_q(f)$, ${L^+_q}={L^+_q}(f)$  and ${L^-_q}={L^-_q}(f)$.

\begin{figure}[tb]
\begin{center}
\psfrag{R}{$R$} %
\psfrag{p}{$(a,c)$} %
\psfrag{q}{$(a,d)$} %
\psfrag{p0}{$(a_0,c_0)$} %
\psfrag{A}{Figure a} %
\psfrag{B}{Figure b} %
\includegraphics[scale=0.35]{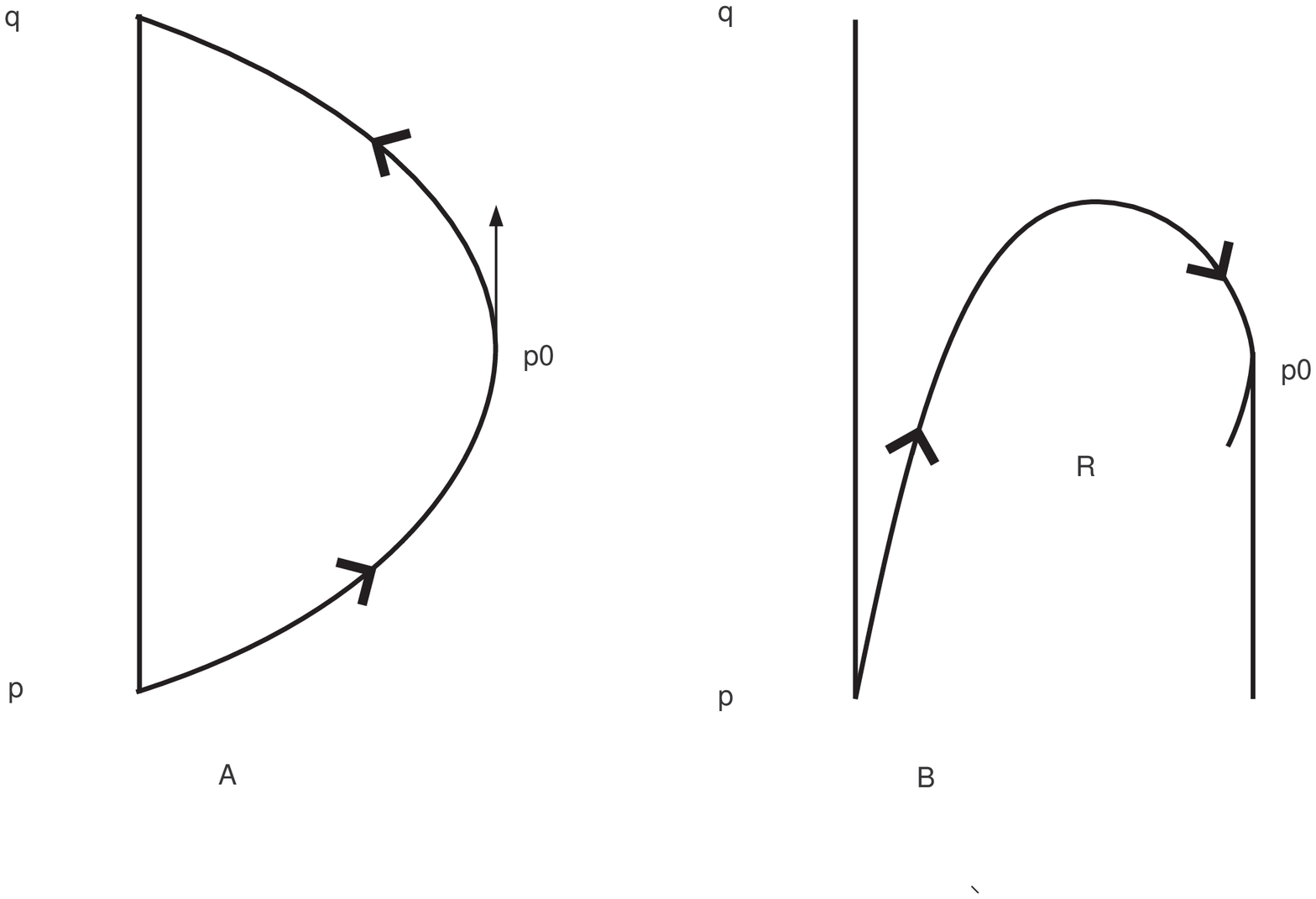}\\
\caption{\footnotesize}\label{fig:disc1}
\end{center}
\end{figure}

\begin{lem}\label{lemma:disk1}
Let $X=(f,g):\mathbb{R}^2\setminus\overline{D}_{\sigma}\to\mathbb{R}^2$ be a differentiable local homeomorphism.
Suppose that the spectrum $\mbox{\rm Spc}(X)\cap[0,+\infty)=\emptyset$ and $p=(a,c)$. Then the intersection of ${L^+_p}$ with the vertical ray $\{(a,y)\in\mathbb{R}^2:y\geq c\}$  is the one point set $\{p\}$.
\end{lem}

\begin{proof}
Assume, by contradiction, that ${L_p^+}\setminus\{p\}$ intersects $\{(a,y)\in\mathbb{R}^2:y\geq{c}\}$.
We take $d>c$ the smallest value such that $q=(a,d)\in{L_p^+}$ (see \mbox{Figure \ref{fig:disc1}a}).
We only consider the case in which the compact arc $[p,q]_f\subset{L_p^+}$ such that $\Pi([p,q]_f)$  equals the interval $[a,a_0]$ with $a<a_0$, where $\Pi(x,y)=x$ (in the other case, $\Pi([p,q]_f)=[a_0,a]$, $a_0<a$ the argument is similar).
Therefore, if we take the vertical segment $[p,q]_a=\{(a,y):c\leq y\leq d\}$ joint to the open disk ${D}(C)$ bounded by the closed curve $C=[p,q]_f\cup[p,q]_a$. We meet two possible cases:

\emph{The first one is that ${D}(C)\cap\overline{D}_{\sigma}=\emptyset$.}
    We select the point  $(a_0,c_0)\in[p,q]_f$ in order that $c_0$ will be the smallest value of the compact set $\Pi^{-1}(a_0)\cap[p,q]_f$.
Let $R$ be the closed region bounded by the union of $\{(a,y):y\leq{c}\},$ $\{(a_{0},y):y\leq {c_{0}}\}$ and $[p,q_0]_f\subset[p,q]_f$ where $q_0=(a_0,c_0)$.
    As $c_0=\inf\{y\in\Pi^{-1}(a_0):(a_0,y)\in[p,q]_f\}$ the compact arc $[p,q]_f$ is tangent to the vertical line $\Pi^{-1}(a_0)$ at the point $(a_0,c_0)$. Thus, $X_f(a_0,c_0)$ is vertical, and so $f_y(a_{0},c_{0})=0$.
This implies that $f_x(a_{0},c_{0})\in\mbox{\rm Spc}(X)$.
    By the assumptions about $\mbox{Spc}(X)$, $f_x(a_{0},c_{0})<0$ which in turn implies that the arc $[q_0,q]_f\subset[p,q]_f$ must enter into $R$ and cannot cross the boundary of $R$ (see \mbox{Figure \ref{fig:disc1}b}). This contradicts the fact that $q=(a,d)\notin{R}$.

\emph{The second case happens when ${D}(C)\cap\overline{D}_{\sigma}\neq\emptyset$.}
As $[p,q]_f$ is not contained in $\overline{D}_{\sigma},$ either $\sigma<a_0$ or $\sigma=a_0$.
If $\sigma<a_0,$ the vertical line $x=\sigma$ meet $[p,q]_f$ in two different points which define a closed curve as in \mbox{Figure \ref{fig:disc1}a} such that it bounds an open disk  disjoint of $\overline{D}_{\sigma}$.
We conclude by using the proof of the first case. If $\sigma=a_0$, we observe the continuous foliation $\mathcal{F}(f)$ in a neighborhood of $[p,q]_f$ and meet two points $\tilde{p}=(\sigma,\tilde{c})$ and $\tilde{q}=(\sigma,\tilde{d})$ with $\tilde{c}<\tilde{d}$, but
$\Pi([\tilde{p},\tilde{q}]_f)\subset[\sigma,+\infty)$.
It satisfies the conditions of the first case.
Therefore the lemma is proved.
\end{proof}

\begin{rmk}\label{rem:disk1}
Lemma~\ref{lemma:disk1} remains  true, if  we consider the negative leaf  $L_q^-$ starting  at $q=(a,d)$ joint to the vertical ray
$\{(a,y)\in\mathbb{R}^2:y\leq d\}$.
\end{rmk}
\begin{figure}[tb]
\begin{center}
\psfrag{p}{$p$} %
\psfrag{q}{$q$} %
\psfrag{a}{$a$} %
\psfrag{b}{$b$} %
\psfrag{A}{$A$} %
\psfrag{B}{$B$} %
\psfrag{D}{$\overline{D}(C)$} %
\psfrag{al}{$\alpha$} %
\psfrag{Aal}{$(\alpha,c)$} %
\psfrag{ma}{$m_{\alpha}$} %
\psfrag{Pa}{$p_{\alpha}$} %
\psfrag{Qa}{$q_{\alpha}$} %
\psfrag{Grad}{$\nabla_f$} %
\psfrag{Xf}{$X_f$} %
\psfrag{Yc}{$[\Pi(p),\Pi(q)]_c$} %
\psfrag{FiA}{Figure a} %
\psfrag{FiB}{Figure b} %
\includegraphics[scale=0.35]{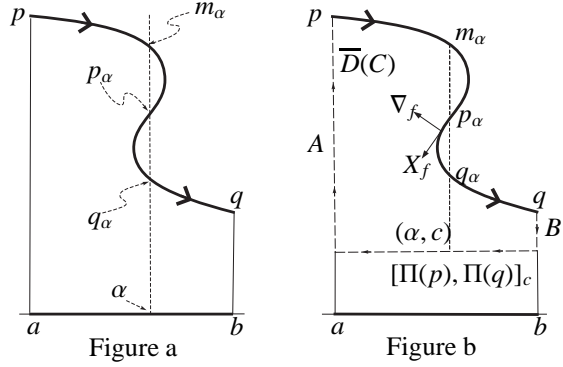}\\
\caption{\footnotesize Here $\big\{\nabla_f(z),X_f(z)\big\}$ a positive basis }\label{fig:disc2}
\end{center}
\end{figure}

\begin{lem}\label{lemma:fluw-positive}
Let $X=(f,g):\mathbb{R}^2\setminus\overline{D}_{\sigma}\to\mathbb{R}^2$ be a map with $\mbox{\rm Spc}(X)\cap[0,+\infty)=\emptyset$.
Consider  $L^+_p\subset\{f=f(p)\}$ and the projection $\Pi(x,y)=x$.
If the oriented compact arc $[p,q]_f\subset{L_p^+}$ and its image $\Pi([p,q]_f)$ is the interval $[\Pi(p),\Pi(q)]\subset(\sigma,+\infty)$ with $\Pi(p)<\Pi(q)$. Then
\[
{\int_{[p,q]_f}}\langle X, \nabla_f\rangle dt \geq f(p) \int_{[p,q]_f}f_x ~ dt+  g(p)\Big(\Pi(q)-\Pi(p)\Big),
\]
where $\langle~,~\rangle$ denotes the usual inner product on the plane, and $f_x$ is the first partial derivative.
\end{lem}

\begin{proof}
For each $\alpha\in \Pi([p,q]_f)=[a,b]$, the vertical line $\Pi^{-1}(\alpha)$ intersects $[p,q]_f$ in a non-empty compact set. So there exist $\hat{y}_{\alpha}=\sup\{y\in\mathbb{R}:(y,\alpha)\in[p,q]_f\}$ and $m_{\alpha}=(\alpha,\hat{y}_{\alpha})$.
We also define  $S\subset(a,b)$ as the set of critical values of  the projection $\Pi(x,y)=x$ restricted to the differentiable arc $[p,q]_f$.
By the Sard's Theorem, presented in \cite[Theorem 3.3]{Gugu} this set $S$ is closed and has zero Lebesgue measure.
For $\alpha\in(a,b)\setminus S$ the set $\Pi^{-1}(\alpha)$ intersects $[p,q]_f$ in at most finitely many points and the complement set $\Big(\Pi^{-1}(\alpha)\cap[p,q]_f\Big)\setminus\{ m_{\alpha}\}$ is either empty or its cardinality is odd.
By Lemma~\ref{lemma:disk1}, the order of these points in the line $\Pi^{-1}(\alpha)$ oriented oppositely to the $y-$axis coincides with that on the oriented arc $[p,q]_f$ (a behavior as in \mbox{Figure \ref{fig:disc1}a} does not exist).
Therefore, for every $\alpha\in(a,b)\setminus S$ the set $\Big(\Pi^{-1}(\alpha)\cap[p,q]_f\Big)\setminus\{ m_{\alpha}\}$ splits into pairs $\{p_{\alpha}=(\alpha,d_{\alpha}),q_{\alpha}=(\alpha,c_{\alpha})\}$ with the following three properties:
    \textbf{(a.1)} $c_{\alpha}<d_{\alpha}$, %
    \textbf{(a.2) }the compact arc $[p_{\alpha},q_{\alpha}]_f$ lies in the semi-plane $\{x\leq\alpha\}$ and it is oriented from $p_{\alpha}$ to $q_{\alpha}$ %
    \textbf{(a.3)} $g(q_{\alpha})>g(p_{\alpha})$.
Notice that the tangent vector of $[p,q]_f$ at $p_{\alpha}$ has a negative $x-$component: $f_y(p_{\alpha})>0$, and the respective tangent vector at $q_{\alpha}$ satisfies $f_y(q_{\alpha})<0$. Similarly, $f_y(m_{\alpha})\leq{0}$ (see \eqref{eq:hamiltonian}).
\begin{description}
    \item[\sc Assertion] \it
    Take $c<\inf\{y:(\alpha,y)\in[p,q]_f\}$ and $[\Pi(p),\Pi(q)]_c=\{(\alpha,c):a\leq\alpha \leq b\}$, an horizontal segment.
    Consider the compact set $\overline{D}(C)\subset\{(x,y):a\leq x\leq b\}$ the boundary of which contain $[p,q]_f\cup[\Pi(p),\Pi(q)]_c$ (see \mbox{Figure \ref{fig:disc2}b}). Suppose that $C$, the boundary of $\overline{D}(C)$ is negatively oriented (clock wise). Then
\[
{\int_{\overline{D}(C)}g_y(x,y)dx\wedge dy}=\int_{[p,q]_f}\langle F, \nabla_f\rangle dt-
{\int_{a}^b g(\alpha,c)d{{\alpha}}},
\]%
where $F(x,y)=\Big(f(x,y)-f(p),g(x,y)\Big)$.
\end{description}

\noindent{\it Proof of Assertion }
We will use the \emph{Green's formulae} given in \cite[Corollary~5.7]{P2} (see also \cite{P1}) with the differentiable map
$G:z\mapsto(0,g(z))$ and the outer normal vector of $C$ denoted by
$z\mapsto\eta(z)$ (unitary). By using that $\mbox{\rm Trace}(DG_z)=g_y(z)$ it follows that
\begin{equation}\label{*}
{\int_{\overline{D}(C)}}g_y(x,y)dx\wedge dy= {\int_{C}}\langle
G,-\eta \rangle ds,
\end{equation}
where $ds$ denotes the arc length element. If $[\Pi(q),\Pi(p)]_c$ denote $[\Pi(p),\Pi(q)]_c$ oriented from $\Pi(q)$ to $\Pi(p)$, then $C=[p,q]_f\cup B\cup[\Pi(q),\Pi(p)]_c\cup A$ with $A$ and $B$ two oriented vertical segments.  Consequently,
\[
{\int_{C}}\langle G,-\eta \rangle ds
={\int_{[p,q]_f}}\langle G,-\eta\rangle ds
+{\int_{B}}\langle G,-\eta \rangle ds
+{\int_{[\Pi(q),\Pi(p)]_c}}\langle G,-\eta \rangle ds
+{\int_{A}}\langle G,-\eta \rangle ds.
\]
In $A\cup B$ the vector $-\eta$ is horizontal
i.e $\eta(z)=(\eta_1(z),0)$ then $G=(0,g)$ implies that $\int_{A}\langle G,-\eta\rangle ds=\int_{B}\langle G,-\eta\rangle ds=0$. Therefore
\begin{equation}\label{**}
{\int_{C}}\langle G,-\eta \rangle ds
={\int_{[p,q]_f}}\langle G,-\eta\rangle ds
+{\int_{[\Pi(q),\Pi(p)]_c}}\langle G,-\eta \rangle ds.
\end{equation}%
In $[p,q]_f$, the outer normal vector es parallel to $-\nabla_f(z)=-\big(f_x(z),f_y(z)\big)$. Then, for all $z\in[p,q]_f$ we obtain that
$-\eta(z)=\frac{\nabla_f(z)}{||\nabla_f(z)||}$ and $G(z)=F(z)$ because $[p,q]_f\subset\{f=f(p)\}$. Thus
\[{\int_{[p,q]_f}}\langle G,-\eta \rangle ds={\int_{[p,q]_f}}\langle F,\nabla_f \rangle dt.\]
Similarly, in $[\Pi(q),\Pi(p)]_c$ we have $-\eta(z)=(0,1)$ thus
\[
{\int_{[\Pi(q),\Pi(p)]_c}}\langle G,-\eta \rangle ds=
{\int_{b}^a g(\alpha,c)d{{\alpha}}}=
-{\int_{a}^b g(\alpha,c)d{{\alpha}}}.
\]
Therefore, \eqref{**} and \eqref{*} prove the assertion. $\Box$
\par
In order to conclude the proof  of this lemma we consider $\overline{D}(C)$ as in the last assertion joint to the construction of its precedent paragraph. Since the complement of $S$ is a total measure set,
\[ %
{\int_{\overline{D}(C)}}g_y(x,y)dx\wedge dy=
{\int_a^b\big(g(m_{\alpha})-g(\alpha,c)\big)d\alpha}+
{\int_{\alpha\not\in{S}}}{\sum_{\Pi(p_{\alpha})=\alpha}\big[g(q_{\alpha})-g(p_{\alpha})\big]}d\alpha.
\] %
Thus, the formulae in the assertion implies that
\begin{equation}\label{***}
\int_{[p,q]_f}\langle F, \nabla_f\rangle dt=
{\int_a^b g(m_{\alpha})d\alpha}+
{\int_{\alpha\not\in{S}}}{\sum_{\Pi(p_{\alpha})=\alpha}\big[g(q_{\alpha})-g(p_{\alpha})\big]}d\alpha.
\end{equation}
But,
\[
\int_{[p,q]_f}\langle F, \nabla_f\rangle dt=\int_{[p,q]_f}\langle X, \nabla_f\rangle dt-f(p)\int_{[p,q]_f}f_x dt
\]
and the property (a.3) of the precedent paragraph  to Assertion~1 shows that  $g(q_{\alpha})-g(p_{\alpha})\geq 0$. Therefore, since $g(m_{\alpha})\geq g(p)$, \eqref{***} implies that
\[
\int_{[p,q]_f}\langle X, \nabla_f\rangle dt-f(p)\int_{[p,q]_f}f_x dt\geq
{\int_a^b g(m_{\alpha})d\alpha}\geq g(p)\big(b-a\big),
\]
and concludes this proof.
\end{proof}

By  applying the methods of the last proof give us the next:

\begin{lem}\label{lem:fluw-negative}
Let $X=(f,g):\mathbb{R}^2\setminus\overline{D}_{\sigma}\to\mathbb{R}^2$ be a map with $\mbox{\rm Spc}(X)\cap[0,+\infty)=\emptyset$.
Consider $L^-_q\subset\{f=f(q)\}$ and $\Pi(x,y)=x$.
If the oriented compact arc $[p,q]_f\subset{L_q^-}$, and $\Pi([p,q]_f)=[\Pi(q),\Pi(p)]\subset(\sigma,+\infty)$ with $\Pi(q)<\Pi(p)$. Then
    \[
    {\int_{[p,q]_f}}\langle X, \nabla_f\rangle dt \geq f(q) \int_{[p,q]_f}f_x dt+
    g(p)\Big(\Pi(p)-\Pi(q)\Big).
    \]
\end{lem}

\begin{proof}
As $a=\Pi(p)<\Pi(q)=b,$ we again consider the null set $S\subset\Pi([p,q]_f)$ given by the critical values of $\Pi$ restricted to $[p,q]_f$ (see \cite{Gugu}).
Similarly, for every $\alpha\in\Pi([p,q]_f)=[a,b]$ we define $\hat{y}_{\alpha}=\sup\{y\in\mathbb{R}:(y,\alpha)\in[p,q]_f\}$ and $m_{\alpha}=(\alpha,\hat{y}_{\alpha})$.
Therefore, Remark~\ref{rem:disk1} shows that for each $\alpha\not\in S$  the finite set $\Big(\Pi^{-1}(\alpha)\cap[p,q]_f\Big)\setminus\{\tilde{m}_{\alpha}\}$ splits into pairs $\{p_{\alpha}=(\alpha,d_{\alpha}),q_{\alpha}=(\alpha,c_{\alpha})\}$ satisfying:
    \textbf{(i)} $c_{\alpha}<d_{\alpha}$, %
    \textbf{(ii)} the oriented arc $[p_{\alpha},q_{\alpha}]_f\subset\{x\geq\alpha\}$, and %
    \textbf{(iii)} $g(q_{\alpha})>g(p_{\alpha})$.
\par
Take $\overline{D}(C)$ the boundary of which is the closed curve $C\subset\Pi^{-1}(b)\cup[p,q]_f\cup\Pi^{-1}(a)\cup[\Pi(p),\Pi(q)]_c$, where $[\Pi(p),\Pi(q)]_c=\{(\alpha,c):a\leq\alpha\leq b\}$ for some $c<\inf\{y:(x,y)\in[p,q]_f\}$.
By using the Green's formulae with the map $z\mapsto(0,g(z))$ and the compact disk $\overline{D}(C)$ we have that
\[
{\int_{\overline{D}(C)}g_y(x,y)dx\wedge dy}=\int_{[p,q]_f}\langle \tilde{F}, \nabla_f\rangle dt-
{\int_{a}^b g(\alpha,c)d{{\alpha}}},
\]%
where $\tilde{F}(x,y)=\Big(f(x,y)-f(q),g(x,y)\Big)$. Since
\[
{\int_{\overline{D}(C)}g_y(x,y)dx\wedge dy}+{\int_{a}^b g(\alpha,c)d{{\alpha}}}
={\int_a^b} g(m_{\alpha}) d\alpha
+{\int_{\alpha\not\in S}} {\sum_{\Pi(p_{\alpha})=\alpha}\big[g(q_{\alpha})-g(p_{\alpha})\big]}d\alpha
\]
and
\[ %
\int_{[p,q]_f}\langle \tilde{F}, \nabla_f\rangle dt
={\int_{[p,q]_f}}\langle X, \nabla_f\rangle dt-f(q) \int_{[p,q]_f}f_x dt,
\] %
the last property \textbf{(iii)} implies
\[
{\int_{[p,q]_f}}\langle X, \nabla_f\rangle dt-f(q) \int_{[p,q]_f}f_x dt
\geq{\int_a^b} g(m_{\alpha}) d\alpha.
\]
This concludes the proof because $g(m_{\alpha})\geq g(p)$ shows ${\int_a^b} g(m_{\alpha}) d\alpha\geq g(p)(b-a)$.
\end{proof}

\section{Hurwitz vector fields}\label{sec:5}

This section concludes with the proof of the main theorem.
The essential goal of the next proposition is to prove the fact that our eigenvalue assumption ensures the non-existence of unbounded half--Reeb components. It is obtained by using the preparatory results of the previous section.
With this fact \mbox{Theorem \ref{thm:main-1}} is just obtained  by applying our previous papers \cite{GR06,GPR06}.
\begin{figure}[tb]
\begin{center}
\psfrag{a}{$a$} %
\psfrag{al}{$\alpha$} %
\psfrag{g0}{$\gamma(0)$} %
\psfrag{g1}{$\gamma(1)$} %
\psfrag{Grad}{$\nabla_f(z)$} %
\psfrag{LM}{${L^+_{\gamma(0)}\subset\{f=0\}}$} %
\psfrag{Lm}{${L^-_{\gamma(1)}\subset\{f=0\}}$} %
\psfrag{mal}{$m_{\alpha}\in\Gamma_{\alpha}^{-}$} %
\psfrag{Sf}{$[\gamma(s),\gamma(\varphi(s))]_f$} %
\psfrag{psal}{$p_{s_{\alpha}}=(a_{s_{\alpha}},d_{s_{\alpha}})\in\{g=g(\gamma(0))\}$} %
\psfrag{qsal}{$q_{s_{\alpha}}=(a_{s_{\alpha}},c_{s_{\alpha}})$} %
\includegraphics[scale=0.5]{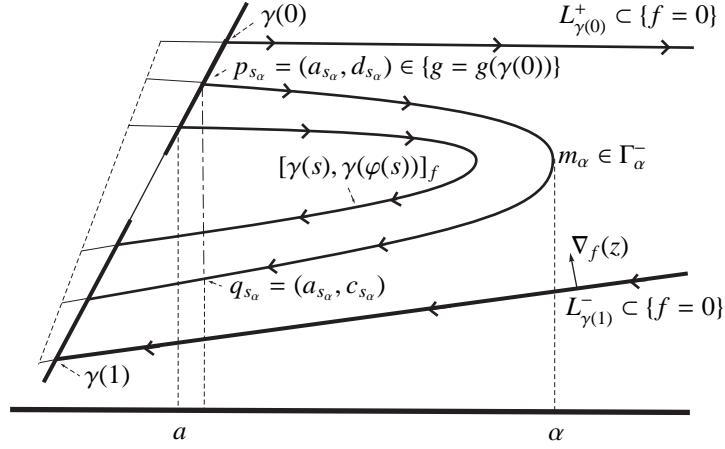}\\
\caption{\footnotesize Half-Reeb componet $\mathcal{A}$, proof of Proposition \ref{prop:main}}\label{fig:HRChurwitz}
\end{center}
\end{figure}
\begin{prop}[Main]\label{prop:main}
Let  $X=(f,g):\mathbb{R}^2\setminus{\overline{D}_\sigma}\to\mathbb{R}^2$ be  a differentiable  map,  where $\sigma>0$ and ${\overline{D}_\sigma}=\{z\in\mathbb{R}^2:||{z}||\leq\sigma\}$.
Suppose that $X$ is Hurwitz: $\mbox{\rm  Spc}(X)\subset\{z\in\mathbb{C}:\Re(z)<0\}$. Then
\begin{enumerate}
\item[(i)] Any half-Reeb component of either $\mathcal{F}(f)$ or $\mathcal{F}(g)$ is a bounded subset of $\mathbb{R}^2$.
  \item[(ii)]There are $s\geq\sigma$ and a globally injective local homeomorphism $\widetilde{X}=(\tilde{f},\tilde{g}):\mathbb{R}^2\to\mathbb{R}^2$ such that $\widetilde{X}$ and $X$ coincide on $\mathbb{R}^2\setminus{\overline{D}_s}$. Moreover, $\mathcal{F}(f)$ and $\mathcal{F}(g)$ have no half-Reeb components.
\end{enumerate}
\end{prop}

\begin{proof}
In the proof of (i), we only consider one case.
Suppose by contradiction that the foliation, given by the level curves has an unbounded half--Reeb component.
By \cite[Proposition 1.5]{FGR1}, there exists a half--Reeb component  of $\mathcal{F}(f)$ the projection of which is  an interval of infinite length. Thus,
\begin{enumerate}
  \item [($a$)]there are $\hat{a}_0>\sigma$ and $\mathcal{A}$, a half--Reeb component  of $\mathcal{F}(f)$ such that $[\hat{a}_0,+\infty)\subset\Pi(\mathcal{A})$, and the vertical line $\Pi^{-1}(\hat{a}_0)$  intersects (transversally) both  non--compact edges of $\mathcal{A}$. Here $\Pi(x,y)=x$.
\end{enumerate}
\par
A half-Reeb component of a fixed foliation contain properly other half-Reeb components, and they are  topologically equivalent.
In this context, the component is stable under perturbations on their compact face as long as the perturbed arc is also a compact face of some component.
Therefore, without lost of generality, we may assume that nearby its endpoints, the compact edge of is made up of arcs of $\mathcal{F}(g)$.
In this way, there exist $0<\tilde{a}<\frac{1}{2}$ and an injective, continuous curve $\gamma:(-\tilde{a},1+{\tilde{a}})\to\mathcal{A}$ such that
\begin{enumerate}
  \item [($b.1$)]  $\gamma([0,1])$ is a compact edge of $\mathcal{A}$ such that both  non--compact edges $L^+_{\gamma(0)}$ and $L^-_{\gamma(1)}$ are contained in a half--plane $\{x\geq \hat{a}\}$, for some $\hat{a}>\sigma$.
  \item[($b.2$)] The images $\gamma\big((-\tilde{a},\tilde{a})\big)$ and $\gamma\big((1-\tilde{a},1+\tilde{a})\big)$  are contained in some leaves of $\mathcal{F}(g)$ such that $\sup\Pi\big(\gamma(1-\tilde{a},1+\tilde{a})\big)<\inf\Pi\big(\gamma(-\tilde{a},\tilde{a})\big)$.
  \item [($b.3$)] For some $0<\delta<\frac{\tilde{a}}{4}$  there exists a orientation reversing injective function
  $\varphi:[-2\delta,2\delta]\to(1-\tilde{a},1+\tilde{a})$ with $\varphi(0)=1$  such that $f\big(\gamma(s)\big)=f\big(\gamma(\varphi(s))\big)$.
  Furthermore, if $s\in(0,2\delta]$ then $\varphi(s)\in(1-\tilde{a},1)$ and there exists an oriented compact arc of trajectory $[\gamma(s),\gamma(\varphi(s))]_f\subset\mathcal{A}$ of $\mathcal{F}(f)$, connecting $\gamma(s)$ with $\gamma(\varphi(s))$.
  \item[($b.4$)]For some $0<\delta<\frac{\tilde{a}}{4}$, small enough  if $s\in(0,\delta]$ and $\gamma(s)=(a_s,d_s)$ then there exists $c_s$ such that $q_s=(a_s,c_s)$ belongs to the open arc $\big(\gamma(s),\gamma(\varphi(s))\big)_f\subset[\gamma(s),\gamma(\varphi(s))]_f$.
   \end{enumerate}
\par
Lemma~\ref{lemma:disk1} implies that
\begin{enumerate}
  \item [($c$)]  for every $s\in(0,\delta]$ as in ($b.4$), $c_s<d_s$.
\end{enumerate}
\par
The eigenvalue condition is invariant under addition of constant vectors to maps, therefore we can assume that
\begin{enumerate}
\item[(d)]$g\big(\gamma(0) \big)>0$, $f=0$ over both  non-compact edges of $\mathcal{A}$
and $f>0$ in the interior of $\mathcal{A}.$
\end{enumerate}
\par
Since $\mathcal{A}$ is the union of an increasing sequence of compact sets bounded by the compact edge and a compact segment of leaf.
Then, from our selection of the compact edge we have that for every $\alpha >\Pi\big(\gamma(0)\big)\geq\inf\{a_s:s\in(0,\delta]\}$, large enough there exists $s_{\alpha}\in(0,\delta]$ as in (b.4) such that the compact arc $[\gamma(s_{\alpha}),\gamma(\varphi(s_{\alpha}))]_f$ projects over $(-\infty,\alpha]$, meeting $\alpha$ (\mbox{Figure \ref{fig:HRChurwitz}}). More precisely,
\[\alpha=\sup\Big\{\Pi(p):p\in[\gamma(s_{\alpha}),\gamma(\varphi(s_{\alpha}))]_f\Big\}.\]
This defines a closed curve $\Gamma^-_{\alpha}$ contained in $\Pi^{-1}(a_{s_{\alpha}})\cup[\gamma(s_{\alpha}),\gamma(\varphi(s_{\alpha}))]_f$.
If $[q_{s_{\alpha}},p_{s_{\alpha}}]_{a_{s_{\alpha}}}\subset\Pi^{-1}(a_{s_{\alpha}})$ is the vertical segment connecting $q_{s_{\alpha}}$, of (b.4)  with $p_{s_{\alpha}}=\gamma(s_{\alpha})$, then this clock wise oriented curve satisfies
\begin{equation}\label{eq:1 main prop}
\Gamma^-_{\alpha}=[q_{s_{\alpha}},p_{s_{\alpha}}]_{a_{s_{\alpha}}}\cup[p_{s_{\alpha}},q_{s_{\alpha}}]_{f}\quad \mbox{and} \quad \Pi(\Gamma^-_{\alpha})=[a_{s_{\alpha}},\alpha],
\end{equation}
where the oriented compact arc $[p_{s_{\alpha}},q_{s_{\alpha}}]_{f}\subset[\gamma(s_{\alpha}),\gamma(\varphi(s_{\alpha}))]_f$.
\par
We select and fix $m_{\alpha}\in\Gamma^-_{\alpha}\cap\Pi^{-1}(\alpha)$, by using \mbox{Lemma \ref{lemma:fluw-positive}} and \mbox{Lemma \ref{lem:fluw-negative}}, respectively we obtain
\[
{\int_{[p_{s_{\alpha}},m_{\alpha}]_f}}\langle X, \nabla_f\rangle dt \geq f(p_{s_{\alpha}}) \int_{[p_{s_{\alpha}},m_{\alpha}]_f}f_x  dt+  g(p_{s_{\alpha}})(x-a_{s_{\alpha}}),
\]
and
\[
{\int_{[m_{\alpha},q_{s_{\alpha}}]_f}}\langle X, \nabla_f\rangle dt \geq f(q_{s_{\alpha}}) \int_{[m_{\alpha},q_{s_{\alpha}}]_f}f_x dt+ g(m_{\alpha})(x-a_{s_{\alpha}}).
\]
Therefore, by adding we conclude
\begin{equation}\label{eq:2 main prop}
{\int_{[p_{s_{\alpha}},q_{s_{\alpha}}]_f}}\langle X, \nabla_f\rangle dt \geq f(p_{s_{\alpha}}) \int_{[p_{s_{\alpha}},q_{s_{\alpha}}]_f}f_x  dt+  g(p_{s_{\alpha}})(x-a_{s_{\alpha}})
\end{equation}
because  $g$ in increasing along $[p_{s_{\alpha}},q_{s_{\alpha}}]_f\subset\{f=f(p_{s_{\alpha}})\}$ and $g(m_{\alpha})\geq g(p_{s_{\alpha}})=g(\gamma(0))> 0$.
\begin{enumerate}
\item[(e.1)] We claim that, the closed curves $\Gamma_{\alpha}^-$, given in \eqref{eq:1 main prop} define the following functions
\[
\alpha\mapsto\left|{\int_{[q_{s_{\alpha}},p_{s_{\alpha}}]_{a_{s_{\alpha}}}}}f dt \right|
\quad \mbox{and} \quad
\alpha\mapsto\left|{\int_{[q_{s_{\alpha}},p_{s_{\alpha}}]_f}}f_x   dt \right|,
\]
they are bounded, when $\alpha$ varies in some interval of infinite length contained on $[\hat{a}_0,\infty)$.
Furthermore,
\[
\lim_{\alpha\to +\infty}f(p_{s_{\alpha}}){\int_{[q_{s_{\alpha}},p_{s_{\alpha}}]_f}}f_x =0.
\]
\end{enumerate}
\par
In fact, by a perturbation in the compact face if it is necessary, it is not difficult to prove that there is some half-Reeb component $\widetilde{\mathcal{A}}$ of $\mathcal{F}(f)$ such that $\widetilde{\mathcal{A}}\supset{\mathcal{A}}$, their boundaries satisfy $\partial \widetilde{\mathcal{A}}\setminus\{\mbox{compact face}\}\supset\partial {\mathcal{A}}\setminus\{\mbox{compact face}\}$ and $\widetilde{\mathcal{A}}\supset[q_{s_{\alpha}},p_{s_{\alpha}}]_{s_{\alpha}}$, for all $s_{\alpha}$.
Since the image $f(\widetilde{\mathcal{A}})\subset f(\mbox{compac face})$, the function $f$ is bounded in the closure of $\widetilde{\mathcal{A}}$, which contain the compact set $\overline{\cup_{s_{\alpha}}[q_{s_{\alpha}},p_{s_{\alpha}}]_{a_{s_{\alpha}}}}$.
Consequently,
\[
\alpha\mapsto\left|{\int_{[q_{s_{\alpha}},p_{s_{\alpha}}]_{a_{s_{\alpha}}}}}f dt \right|
\]
is bounded, in some interval of infinite length contained on $[\hat{a}_0,\infty)$.
In order to prove the second part, we apply the Green's formulae to the map $z\mapsto(1,0)$ in the compact disk $\overline{D}(\Gamma_{\alpha}^-)$. Since the trace is cero, we obtain
\[
\left|{\int_{[q_{s_{\alpha}},p_{s_{\alpha}}]_f}}f_x   dt\right|
=\mbox{arc length of } [q_{s_{\alpha}},p_{s_{\alpha}}]_{a_{s_x}}.
\]
By compactness we conclude. The last part is directly obtained by using \eqref{eq:1 main prop} and $p_{s_{\alpha}}=\gamma(s_{\alpha})$  because the continuity of the foliation $\mathcal{F}(f)$ and (d) imply that
\[
\lim_{\alpha\to +\infty}p_{s_{\alpha}}=\gamma(0)\in\{f=0\}.
\]
Therefore (e.1) holds.
\begin{enumerate}
\item[(e.2)] We claim that,
\[
\mbox{if}\quad\alpha\rightarrow+\infty
\quad\mbox{then}\quad
{\oint_{\Gamma_{\alpha}^-}}\langle X, \eta^i_{\alpha}\rangle dt\rightarrow+\infty,
\]
where $-\eta^i_{\alpha}$ is a outer normal vector of the close wise oriented curve $\Gamma_{\alpha}^-$.
\end{enumerate}
\par
In the compact arc $[p_{s_{\alpha}},q_{s_{\alpha}}]_f\subset\Gamma_{\alpha}^-$, the vector $\eta^i_{\alpha}$ is parallel to $\nabla_f$.
Thus \eqref{eq:2 main prop} and (e.1) imply that
\[
{\oint_{\Gamma_{\alpha}^-}}\langle X, \eta^i_{\alpha}\rangle dt\geq
A+g(p_{s_{\alpha}})(\alpha-a),
\]
where the constant $A$ is independent of $\alpha$ and $a=\min\{a_s:s\in(0,\delta]\}$.
Since (d) and (b.2) imply that $g(p_{s_{\alpha}})=g(\gamma(0))>0$, it is not difficult to  obtain (e.2).
\par
By (e.2) we select some $\tilde{\alpha}>a$ such that  $\Gamma_{\alpha}^{-}$ satisfies
$\displaystyle{\oint_{\Gamma_{\tilde{\alpha}}^{-}}}\langle X, \eta^i_{\alpha}\rangle dt>0$.
By using the Green's formulae with the map $X=(f,g)$, the assumptions over the eigenvalues i.e $0>\mbox{\rm
Trace}(DX_z),$ imply that
\[
0>\displaystyle{\oint_{\Gamma_{\tilde{\alpha}}^{-}}}\langle X, \eta^i_{\alpha}\rangle dt>0.
\]
This contradiction concludes the proof of part (i).
\par
In order to obtain (ii), we apply \mbox{Theorem \ref{teo:1}} because $X$ satisfies its conditions.
This gives the existence of the pair $(\tilde{X},s)$ with $s\geq\sigma$ and $\tilde{X}=(\tilde{f},\tilde{g}):\mathbb{R}^2\to\mathbb{R}^2$ a globally injective local homeomorphism such that $\widetilde{X}$ and $X$ coincide on $\mathbb{R}^2\setminus{\overline{D}_s}$.
Furthermore, the last property in (ii) is obtained as a direct application of \mbox{Corollary \ref{cor:ext with no HRC}}. Therefore, this proposition holds.
\end{proof}

Now we prove our main result

\subsection{Proof of Theorem \ref{thm:main-1}}
By \mbox{Proposition \ref{prop:main}}, there exists a globally injective local homeomorphism $\widetilde{X}:\mathbb{R}^2\to\mathbb{R}^2$ such that $\widetilde{X}$ and $X$ coincide on some  $\mathbb{R}^2\setminus{\overline{D}_{s_1}}$, with $s_1\geq\sigma$.
In particular, the restriction  $X|:{{\mathbb{R}^2\setminus \overline{D}_{s_1}}}\to \mathbb{R}^2$ is injective.
In order to shown the existence of the differentiable extension, consider $v=-\widetilde{X}(0)$ joint to the globally injective  map $\widetilde{X}+v$.
In this context, we can apply the arguments of \cite[Theorem 11]{GPR06}. Thus there are $\tilde{s}>s_1\geq\sigma$ and a global differentiable  vector field
$Y:\mathbb{R}^2\to\mathbb{R}^2$ such that
\begin{enumerate}
    \item[(a.1)] ${\mathbb{R}^2\setminus\overline{D}_{\tilde{s}}}\ni z\mapsto Y(z)=(X+v)(z)$  is also injective and $Y(0)=0$.
    \item[(a.2)] The map $z\mapsto\mbox{\rm Trace}(DY_z)$ is Lebesgue almost--integrable in whole $\mathbb{R}^2$
    (\cite[\mbox{Lemma 7}]{GPR06}).
    \item[(a.3)] The index $\mathcal{I}(X+v)$ is a well-defined number in $[-\infty,+\infty)$
    (\cite[\mbox{Corollary 13}]{GPR06}).
\end{enumerate}
Thus, there exist the index of $X$ at infinity, $\mathcal{I}(X)=\mathcal{I}(X+v)$. Therefore, $\widehat{X}=Y-v$ is the global differentiable extension of $X|:{{\mathbb{R}^2\setminus \overline{D}_{\tilde{s}}}}\to \mathbb{R}^2$ and the pair $(\widehat{X},\tilde{s})$ satisfies the definition of the index of $X$ at infinity. This concludes the proof of (i).
\par
To prove the first part of (ii) we refer the reader to \cite[Lemma 3.3]{FGR1}. Furthermore, since $\mbox{Trace}(D(X+\tilde{v}))=\mbox{Trace}(DX)<0$, for every constant vector $\tilde{v}\in\mathbb{R}^2$ we obtain that
\begin{enumerate}
\item[(b.1)] Given a constant $\tilde{v}\in\mathbb{R}^2$, the vector field $X+\tilde{v}$ generates a positive semi-flow on $\mathbb{R}^2\setminus D_{\sigma}$.
\end{enumerate}
\par
An immediate consequence of (i) is that: if $X$ is Hurwitz, then outside a larger disk both $X$ and $Y$ have no rest points.
In addition, by (a.1) the Hurwitz vector field $Y$ has no periodic trajectory $\gamma$ with $D(\gamma)$ contained in $\mathbb{R}^2\setminus\overline{D}_\sigma$.
As $\mbox{Trace}(DY)=\mbox{Trace}(DX)<0$ by Green's Formulae $Y$ admits at most one periodic trajectory, say $\gamma,$ such that $\overline{D}(\gamma)\supset\overline{D}_\sigma$. Consequently
\begin{enumerate}
\item[(b.2)]
    There exit $s>\tilde{s}$ such that $Y$ satisfies (a.1), (a.2) and $\mathbb{R}^2\setminus D_s$ is free of rest points and periodic trajectories of $Y$.
\end{enumerate}
Under these conditions (b.1) and  \cite[Theorem 26]{GPR06} imply that:
\begin{enumerate}
\item[(b.3)]
    For every $r\geq s$ there exist a closed curve $C_r$ transversal to $Y$ contained in the regular set $\mathbb{R}^2\setminus \overline{D}_r$. In particular, $D(C_r)$ contains $D_r$ and $C_r$ has transversal contact to each small local integral curve of $Y$ at any $p\in C_r$.
\end{enumerate}
Moreover, \cite[Theorem 28]{GPR06} shown that:
\begin{enumerate}
\item[(b.4)] The point at infinity of the Riemann Sphere $\mathbb{R}^2\cup\{\infty\}$ is either an attractor or a repellor of $X+v:\mathbb{R}^2\setminus\overline{D}_s\to\mathbb{R}^2$. More specifically, if $\mathcal{I}(X)<0$ (respectively $\mathcal{I}(X)\geq0$), then $\infty$ is a repellor (respectively  an attractor) of the vector field $X+v$.
\end{enumerate}
Therefore, (ii) holds and concludes the proof of {Theorem \ref{thm:main-1}}. $\Box$


\end{document}